\documentclass[12pt]{article}
\usepackage{amssymb,amsbsy,amsmath,amsfonts,amssymb,amscd,amsthm}
\usepackage{latexsym}
\def\k{\bf k}

\def\phi{\varphi}

\def\a{\alpha}

\def\Z{{\mathbb Z}}

\newtheorem{theorem}{Theorem}[section]
\newtheorem{lemma}[theorem]{Lemma}
\newtheorem{corollary}[theorem]{Corollary}

\newtheorem{definition}[theorem]{Definition}

\begin{document}

\title{Non-trivial stably free modules over crossed products }

\author{Natalia~K.~Iyudu$^1$
and Robert~Wisbauer$^2$ }

\date{}

\maketitle

{\small\tt\noindent \llap{$^1\,$} Department of Pure Mathematics,
Queen's University Belfast, UK \\

{\small\tt\noindent\llap{$^2\,$} Mathematisches Institut, Heinrich
Heine Universit\"at, \hfill\break\noindent 40225 D\"usseldorf,
Germany

}

\medskip

\bf \rm e-mails:\
n.iyudu@qub.ac.uk,\\
\hglue\parindent\phantom{{\bf e-mails:\
}}wisbauer@math.uni-duesseldorf.de

\begin{abstract}

We consider the class of crossed products of noetherian domains with
universal enveloping algebras of Lie algebras. For algebras from
this class we give a sufficient condition for the existence of
projective non-free modules. This class includes Weyl algebras and
universal envelopings of Lie algebras, for which this question,
known as noncommutative Serre's problem, was extensively studied
before. It turns out that the method of lifting of non-trivial
stably free modules from simple Ore extensions can be applied to
crossed products after an appropriate choice of filtration. The
motivating examples of crossed products are provided by the class of
RIT algebras, originating in non-equilibrium physics.

\rm\bigskip \noindent {\bf Keywords:} \ Serre's problem, stably free
modules, Ore extension, crossed product with universal enveloping of
Lie algebra, RIT (relativistic internal time) algebras, faithfully
flat modules, type of the element, strongly completely prime
subalgebra, graded domain, ordered-like semigroup
\end{abstract}

\section{Introduction}

In  \cite{serre} J.-P.Serre posed  the question on whether any
finitely generated projective module over the ring of commutative
polynomials ${\k}[x_1,\dots,x_n]$ over a field $\k$ is free. It was
stated there in  geometrical language: whether any locally trivial
vector bundle over an affine space ${\mathbb A}^n_{\k}$ is a trivial
bundle. After almost twenty years of attempts A.~Suslin \cite{sus}
and D.~Quillen \cite{qui} independently (and using different
methods) obtained an affirmative answer to  the Serre question (see
also \cite{lam} for a detailed study of the techniques involved).

Later on, this question was investigated for various classes of
non-commutative rings.  A report on some of this work also can be
found in Lam's book \cite{lam}, ch.VII.8. To describe briefly what
has been done let us remind some definitions.

 A finitely generated left $A$-module $M$ is
called {\it stably free} if $M\oplus A^n=A^m$ for some nonnegative
integers $n$ and $m$, clearly, it is projective then. A module
which is stably free but not free will be called {\it non-trivial
stably free}. We will need the {\it rank}
of the stably free module $M$ which is defined as ${\rm rk} M =
m-n$. This definition obviously only makes sense   if $A$ is an
IBN-ring, for example, we may consider noetherian rings. We will
suppose throughout the paper that all rings are IBN.

 Situation in non-commutative case turned out to be more
involved: there were constructed counterexamples, i.e. stably free
non-free modules in several classes of non-commutative algebras. By
saying that there exists a counterexample in a certain class of
algebras, we mean that any algebra from this class allows a finitely
generated projective but non-free module.

For example, stably free non-free ideals were constructed by Webber
\cite{web}  in any Weyl algebra $A_n$. Another counterexample was
constructed by Ojanguren and Sridharan \cite{os} in rings of
polynomials on two variables over a division ring (which is not a
field). In group algebras, non-free projective modules were
constructed by Dunwoody and Berridge \cite{DuBe} for torsion free
groups and by Artamonov  \cite{Art81} for solvable groups. Examples
of this type in enveloping algebras of non-Abelian finite
dimensional Lie algebras were provided by  Artamonov
%\cite{Art85}
and by Stafford,  in  \cite{staff}  a unified way for producing
non-trivial stably free right ideals was given, which virtually
covers above cases.

In this note we consider the class of crossed products of noetherian
domains with a universal enveloping algebra of a Lie algebra, which
subsumes most of classes mentioned above, and provide a sufficient
condition for the existence of stably free non-free modules in this
wider class. More precisely, we show in theorem \ref{t3} that stably
free non-free modules can be lifted from any subalgebra of the
crossed product $ A \star U{\cal G}$ which is a simple differential
Ore extension $A[g,\delta]$, $g \in {\cal G}$, $\delta \in Der A$
($\delta=\delta_{\bar g}$ is a derivation, involved in the given
crossed product, associated to the element $g \in {\cal G}$).

As an element of our tool we prove the following (theorem
\ref{ft3}): if $A$ is a domain, endowed with a filtration by a
well-ordered semigroup, such that $A_0$ is a faithfully flat
$A$-module, then any non-trivial stably free ideal in $A_0$ could be
lifted to a nontrivial stably free ideal in $A$. Graded version of
this fact (theorem \ref{gt3}) we prove for a graded domain $A$,
graded by {\it order-like semigroup}. This is a wider class of
semigroups, which however captures most essential properties of
well-ordered semigroups.

Let us emphasize   that all examples of non-trivial stably free
modules mentioned above, and just about all known examples in the
noncommutative case, are modules of rank one. Over commutative
rings, there are examples of higher rank and these are typical. An
example of module of minimal rank (over commutative rings) is a
module of rank two over the ring ${\mathbb R}[x,y,z]/x^2+y^2+z^2=1$.
The way to ensure that the right unimodular row $(x,y,z)$ gives rise
to the non-trivial stably free module has a geometrical flavour
(using the theorem of the "hedgehog brushing" on a 2-sphere) and
does not explain much in the line of techniques we study here.

For the class of Weyl algebras $A_n(k)$ it was proved by Stafford
\cite{StW} that all stably free modules of rank two and bigger are
free. In the proof, of course, simplicity of $A_n(k)$ plays a
crucial role.  This result was generalized to some crossed products
of simple rings with supersolvable groups by Jaikin-Zapirain in
\cite{J-Z}. Another cases where positive result holds one can find
in \cite{artc}.

Examples  of crossed products with the universal enveloping algebra
we consider are provided by the class of RIT (relativistic internal
time) algebras, which we have been studying in \cite{wia},
\cite{ia}. This class originated in non-equilibrium physics
\cite{am}, \cite{a1} and we consider it here in the general setting
of crossed products.

%As in general, crossed products of this type are common in
%constructions originating in physics, the class of RIT algebras
%appeared in \cite{am}, \cite{a1} in relation to certain theories of
%non-equilibrium physics.

\section{Choice of subalgebras}\label{int}

 We start with the definition of a crossed product with a
universal enveloping algebra.

{\bf Definition.} Let $A$ be a (noncommutative) $\k$-algebra,
and ${\cal G}$ a Lie algebra  with a basis $\{ g_i | i \in I\}$ over
$\k$. Then a $\k$-algebra $B$ containing $A$ is called a {\it
crossed product} provided there is an embedding ${\cal G}
\hookrightarrow B: g \mapsto \bar g$ of linear spaces, which
satisfies:

\vskip3mm

(1)  $\bar g  r - r \bar g \in A$ for any $g \in G, r\in A$

and $\bar g r - r \bar g=\delta_g(r)$ is a derivation on $A$.

\vskip3mm

(2) $ \bar g  \bar h - \bar h \bar g = \overline{[g,h]} $ mod $A$
for all $g,h \in {\cal G}\,\,\,$

\vskip3mm

(3) $B$ is a free (right)  $A$-module with the commutative monomials
$\bar g_1^{j_1}...\bar g_m^{j_m}$ on $\{ g_i | i \in I\}$ as a basis

\vskip5mm

We denote the crossed product algebra $B$ by $ A \star U{\cal G}$.

 It is known from \cite{staff} (see also \cite{McR}) that
non-trivial stably free ideals do exist in  simple differential Ore
extensions of noetherian domains, which satisfy some additional
condition.

In this section we suggest how to choose appropriate subalgebras in
a crossed product $ A \star U{\cal G}$ in such a way that their
non-trivial stably free modules can be lifted to non-trivial
stably free modules over the whole crossed product. Namely, we take
subalgebras isomorphic to a simple Ore extension of the initial
algebra $A$, thus the idea is to use subalgebras in the "intersection"
of the crossed product components.

We prove several properties of these subalgebras in order to
prepare the tool which allows us to lift non-trivial stably free
modules from these subalgebras to the whole crossed product.

Directly from the definitions, it can be seen that $A
\star U{\cal G}$, where ${\cal G}=\{g\}$ is a one-dimensional Lie
algebra,  is isomorphic to the simple differential Ore
extension $A[x, \delta]$, where $\delta=\delta_g$ is the derivation
related to $g$, which was defined above as: $\,\, \delta_g(r)=g \bar
r- \bar r g,$    for $ r \in A$.

Using the defining relations (1) and (2) in the crossed product one can
easily see that $A \star U{\cal G}_1$, where ${\cal G}_1=\{g\}$ is a
Lie algebra generated by any single element $g \in \cal G$, is a
subalgebra in $A \star U{\cal G}$. Indeed, the free basis of the
$A$-module $A \star U{\cal G}_1$ according to (3) consists of
elements $\bar g^i, \, i=0 \leq i < \infty $. Any product of two
elements of the shape $a \bar g^i, \, a \in A, \, g \in \cal G$ is a
linear combination of elements of the same shape after applying (1).

We denote by $A_1$ the subalgebra $A \star U{\cal G}_1$ in $B=A
\star U{\cal G}$ and will consider $B$ as a  left
$A_1$-module writing ${}_{A_1} B$.

\section{Faithful flatness of ${}_{A_1} B$}\label{ff}

Using the above notations we will prove two crucial properties of
subalgebras of our choice which allow to lift nontrivial stably free
modules from $A_1=A \star U{\cal G}_1=A[g,\delta]$ to the whole
crossed product $A \star U{\cal G}$.

Starting from here we suppose that $A$ is Noetherian domain.

The first property we need is the faithfully
flatness of $B$ as a left $A_1$ module. We will prove that in our
situation even a stronger condition holds, namely

\begin{lemma}\label{l1} The left $A_1$-module  $ {}_{A_1} B$ is free.
\end{lemma}

\begin{proof}
 By definition of a crossed product,  $B=A
\star U{\cal G}$ and $A_1=A \star U{\cal G}_1$ are free left
$A$-modules with bases $V=\{ g_1^{i_1} \dots g_n^{i_n}   \}$ and
$W=\{ g_1^{i} \}$ respectively.

 We prove that $ {}_{A_1} B$ is generated by the set
$\Omega = g_2^{j_2} \dots g_n^{j_n}$ and this set forms a free basis
of this $A_1$-module.

The first part of the statement, saying that $\Omega$ is  generating system is obvious.
To show that this is a free basis it is enough to check
%For the second one is equivalent to show
that if $ \Sigma a_ib_i = 0 $ in $B$, for $ a_i \in A_1, b_i
\in \Omega$, then all $a_i=0$, since $A_1$ has no zero divisors. The
equivalence follows from our condition that $A$ is a domain.
%{\bf where is this requested ?}
To ensure this let us first write elements $a_i \in A_1$ from the
sum $ \Sigma a_ib_i $ above as follows:
%{\bf what does this mean ?}
$a_i=a_i(g_1)=\Sigma \alpha_k^{(i)} g_1^k$, with $ \alpha_k^{(i)}
\in A $.

 Now fulfill the multiplication in the above sum and gather terms
near each $v_i \in V$. We get $ \Sigma  \beta_j v_j = 0 $, where
$\beta_j=\Sigma_{i,k: \, v_j=g_1^{k}v_i} \,\, \alpha_k^{(i)}$. From
this it follows that $\beta_i=0$, since $V$ is a free generating
system for ${}_A B$.
Since $V$
is a free basis for $A$-module  $B$, given the fixed numbers $i$ and
$j$, there is only one $k$ such that $v_j=g_1^{k}v_i$. Hence the set
$\{ \beta_j \}$ just coincides with the set $\{ \alpha_k^{(i)} \}$.
So, together with all $\beta_i=0$ we have all $\alpha_k^{(i)}=0$, and
hence $a_i=0$ for all $i$.
\end{proof}
\medskip

\section{Strongly completely prime subalgebras}\label{scp}

Before we start the discussion of the second main lemma we should
introduce the notion of {\it strongly completely prime subalgebra},
or {\it s.c.p.-subalgebra} for short.

Let us consider the following two properties of subalgebra.

{\bf Definition 1.} We say that a {\it subalgebra} $A_1$  is {\it
completely prime} in  $A$ if for any two non-zero elements $a$ and
$b$ from $A$,
 $ab \in A_1$ implies  $a \in A_1$ or $b \in A_1$.

{\bf Definition 2.} We say that a {\it subalgebra} $A_1$  is {\it
strongly completely prime ($s.c.p.$)} in  $A$ if
 for any two non-zero elements $a$ and
$b$ from $A$, $ab \in A_1$ implies  $a \in A_1$ and $b \in A_1$.

In case $A_1$ is an ideal in $A$, the first definition just coincide
with the definition of a completely prime ideal, e.i. an ideal such
that the quotient is a domain. (this explains the origin of our
term).

The second definition degenerate in case $A_1$ is an ideal. Indeed,
suppose $A_1 \triangleleft A,$ $\, A_1 \neq \{0\}$ and $A \backslash
A_1 \neq \emptyset $. Then take an element $b \in A \backslash A_1$,
since $A_1 $ is a (right) ideal, for arbitrary nonzero element $a
\in A_1$ we have $ab \in A_1$ and the property from the definition 2
doesn't hold. If $A_1= \emptyset$ then formally property of being
$s.c.p.$ always holds in a domain.

So, property of being $s.c.p.$ is clearly a feature of subalgebras
and should be considered only in this case (rather then for ideals).

Let us discuss now the notion of type of an element in the crossed
product.

Firstly, we associate with any product (monomial) $w=ag_{i_1}\dots
g_{i_n}$, $i_k\in I$, $a\in A$ in the crossed product $B=A*{\cal
U}({\cal G})$ its type on variables $g_{i_1}\dots g_{i_n}$. By
definition  {\it the type'} $t(w)$ of the element $w$ is a tuple of
nonnegative integers $(j_1,\dots,j_r)$, where $j_k$ is the number of
variables $g_k$ in the monomial $w$ for any $k\in I$. (In case
$j_l=0$ for all $l>r$, $j_l=0$, we just omit zero terms in the
sequence $j_1,\dots,j_r,\dots$ starting from $j_{r+1}$ to get
$t(w)$). One can also consider the type of a product on any subset
of variables $\{g_{i_k}, \,\,\, i_k \in I' \subset I\}$.
%{\bf is the notion of "type" standard ?}

Let us fix an order on monomials $w=ag_{i_1}\dots
g_{i_n}\in B$ using the degree lexicographical ordering on
commutative words $t(w)$. Namely, we say that $w>w'$ for
$w=ag_{i_1}\dots g_{i_n}$ and $w'=a'g_{i'_1}\dots g_{i'_{n}}$ if
$t(w)\mathop{>}\limits_{dl}t(w')$. The latter means that if
$t(w)=(j_1,\dots,j_r)$ and  $t(w')=(j'_1,\dots,j'_s)$, then either
$r>s$ or $r=s$ and $j_t>j'_t$ for some $t$, such that $j_l=j'_l$ for
all $l\leq t$ .

We  can define a {\it normal form} (with respect to $g_i$, $i\in I$)
of an element in $B=A*{\cal U}({\cal G})$. We say that an element
$f=\sum a_{\overline{i}}g_{i_1}\dots g_{i_n}$ is {\it in the normal
form} if $i_1\leq i_2\leq {\dots}\leq i_n$, that is, all monomials
have the form $g_1^{j_1}\dots g_r^{j_r}$ with coefficients from $A$.
It is clear from the relations in the definition of crossed product
that any element from $B$ can be presented in a normal form, since
these relations allow to commute $r \in A$ with $g \in {\cal G}$ and
elements from ${\cal G}$ between each other. In both cases we might
get new terms, which have a lower degree in $g_i$, $i \in I$. Since
there is no infinite chain of words in $g_i$ of strictly increasing
degree, in certain step we will get an element equal to $w$ in a
normal form. This element we will call a {\it normal form of } $w
\in B=A*{\cal U}({\cal G})$ and denote it by ${\cal N}(w) $.
Property (3) in the definition of the crossed product (PBW -
property) ensures that the normal form of the element in $B$ is
unique.

This allows us to introduce the notion
of the type of an element $b\in B$.

{\bf Definition.} By the { \it type} of an arbitrary element $b\in
B$, we call the type of the highest monomial in the normal form of
$b$.

Having in hands the notion of the type of an element in $B=A \star
U{\cal G}$ we actually have a natural filtration on $B$. Namely $B=
\mathop{\cup}\limits_{\bar i \in \Sigma_n} B_{\bar i}$, where
$\Sigma_n$ is a semigroup of tuples $(i_1,...,i_n)$ with the
componentwise operation and $B_{\bar i}$ is a linear span over $\k$
of elements of the type $(i_1,...,i_n)$, in particular, $B_{ 0}=A$.

The existence of such a filtration force us to develop a general
machinery for the graded and filtered case and then apply it to the
situation of crossed products, using however a filtration different
from the above.

\section{
Semigroup graded and filtered case \label{gr}}\rm

For this section we break our agreement that $A$ is a domain, in
some statements here we will ask only for $A$ being a graded domain
(that is, there are no zero divisors among homogeneous elements with
respect to a given grading).

The main theorem in the graded setting will have the form.

\begin{theorem}\label{gt2} Let $A=\bigoplus_{j\in\Z_+}A_j$ be a
$\Z_+$-graded domain, where $A$ is a flat $A_0$-module. Then any
stably free non-free module over $A_0$ can be lifted to $A$.
\end{theorem}

This theorem can be further generalized in a sense that one can
consider gradings more general than $\Z_+$-gradings. We shall prove
a theorem in that bigger generality, so Theorem~\ref{gt2} will
follow from Theorem~\ref{gt3}.

\begin{definition}\label{df1} We call a semigroup $(G,+)$ {\it ordered-like}
if it has no invertible elements except 0 and for any two finite
subsets $S_1,S_2$ of $G$ such that $S=S_1+S_2\neq\{0\}$, there
exists $c\in S$ with
$$
\nu(c)=|\{(a,b):a\in S_1,\ b\in S_2,\ a+b=c\}|=1.
$$
\end{definition}

Most common example of a semigroup with such a property is a
well-ordered semigroup, where there exists a linear order compatible
with an operation: $a < b \Longrightarrow a+c < b+c $. In this case
as an element $c \in S $ with unique presentation as a sum of
elements from $S_1$ and $S_2$ ($|\nu(c)|=1$) will serve a sum of
maximal elements of $S_1$ and $S_2$. But there are other examples
where this property doesn't come from well-ordering.

\begin{theorem}\label{gt3} Let $A=\bigoplus_{\sigma\in G}A_\sigma$ be
a domain graded by an ordered-like semigroup $G$, where $A$ is
faithfully flat as a left $A_0$-module. Then any stably free
non-free right ideal in $A_0$ can be lifted to $A$.
\end{theorem}

The proof is based on the following lemmas.

\begin{lemma}\label{le1} Let $A=\bigoplus_{\sigma\in G}A_\sigma$ be
a  graded domain, with $G$ being ordered-like semigroup. Then $A_0$
is a completely prime subalgebra of $A$.
\end{lemma}

\begin{proof} To prove that $A_0$ is  completely prime it is enough to
ensure that for any two elements $a,b\in A$, from $a,b\notin A_0$ it
follows that $ab\notin A_0$. For any element $a$ of $A$ let us
denote by $S(a)$ the subset $S(a)=\{\sigma\in G:a_\sigma\neq 0\}$ of
the semigroup $G$, where $a=\sum\limits_{\sigma\in G}a_\sigma$,
$a_\sigma\in A_\sigma$ is the graded decomposition of $a$. Clearly
$a=\sum\limits_{\sigma\in S(a)}a_\sigma$. Since $a,b\notin A_0$, the
sets $S(a)$ and $S(b)$ contain non-zero elements. Since $G$ has no
non-zero invertible elements, we have $S=S(a)+S(b)\neq\{0\}$. Taking
into account that $G$ is ordered-like, we can find $\gamma\in
S\setminus\{0\}$ such that $ \nu(\gamma)=|\{(\sigma,\tau):\sigma\in
S(a),\ \tau\in S(b),\ \sigma+\tau=\gamma\}|=1$. For the $\gamma$-th
graded component of $ab$ we will have $(ab)_\gamma=a_\sigma b_\tau$.
Since $a_\sigma\neq 0$, $b_\tau\neq 0$ and $A$ is a graded domain,
we get $(ab)_\gamma\neq 0$ for $\gamma\neq 0$. So, $ab\notin A_0$.
\end{proof}

The following fact is true for grading by any semigroup, not
necessarily with the ordered-like property.

\begin{lemma}\label{le2} Let $A=\bigoplus_{\sigma\in G}A_\sigma$ be
a domain graded by an arbitrary abelian semigroup $G$. Then the
property of $A_0$ to be completely prime implies the property of
$A_0$ to be strongly completely prime.
\end{lemma}

\begin{proof} To ensure this we should show that if $a,b\in A\setminus
\{0\}$ and $ab\in A_0$ implies $a\in A_0$, then we also have $b\in
A_0$.

Indeed, let $b=\sum\limits_{g\in S(b)}b_g$ be the graded
decomposition of $b$. Then $ab=\sum\limits_{g\in S(b)}ab_g$. Here
$(ab)_g=ab_g\in A_0A_g\subseteq A_g$. On the other hand, $ab\in A_0$
and therefore $(ab)_g=ab_g=0$ for any $g\neq 0$.  Since $a\neq 0$
and $A$ is a graded domain, this implies that $b_g=0$ for any $g\neq
0$. That is, $b\in A_0$.
\end{proof}

As a corollary of Lemmas~\ref{le1} and~\ref{le2} we have that the
subalgebra $A_0$ of $A$ is strongly completely prime. Using this we
can proceed with the proof of the theorem\ref{gt3} by analogy with
\cite{staff}.

\begin{proof}(of the Theorem \ref{gt3})

Let $K$ be a  nontrivial stably free right ideal in $A_0$. We will
show that the induced ideal $KA=K \otimes_{A_0} A$ is also stably
free but not free.

Since $A$ is flat as $A_{0}$-module, $KA=K \otimes_{A_0} A$ and is
also projective as $A$-module, hence stably free. The essential part
is to prove that it is not free. We have $KA \oplus A = A \oplus A$,
thus we have to show that $KA$  is not cyclic.

Suppose this is not the case, i.e. $KA= yA$ for some $0 \neq y \in
A$. (In case $y=0$ we will have a contradiction immediately: $A=yA
\oplus A = A \oplus A$ and this contradicts with the condition we
suppose to holds throughout   the paper that all rings have the IBN
property).

Since $K$ sitting inside $KA$ ($A$ has a unit and all modules are
unital) $K \subset  KA  = yA$, we can take a nonzero element $p \in
K$, which is $p=yb$ for some nonzero $b \in A$. But $K$ is an ideal
in $A_{0}$ and we can use primeness of $A_{0}$: if the element $p
\in A_{0}$ and $p=yb$ for nonzero $y,b \in A$ then it should imply
$y \in A_{0}$.

Now, for any right ideals $I \varsubsetneqq J \triangleleft_r
A_{0}$, due to faithfully flatness of $A$ we have $IK \varsubsetneqq
JK \triangleleft_r A$.

Suppose  that the following inclusion of right ideals in $A_{0}$
holds: $K \varsubsetneqq KA \cap A_{0}$. Then applying the above
observation we
 get $KA \varsubsetneqq (KA \cap A_{0})A$. But in fact $(KA \cap
A_{0})A=KA$: $K \subset KA \cap A_{0}$ implies $KA \subset (KA \cap
A_{0})A$, while $KAA \cup A_{0}A \subset KA \cap A = KA$, hence
$K=KA \cap A_{0}$. This contradiction shows  that $K=KA \cap A_{0}$.

Now we have $K=KA \cap A_{0}=yA \cap A_{0} \subset yA_{0}$, due to
$yb \in yA$ belongs to $A_{0}$ we again use the fact that $A_{0}$ is
a prime subalgebra, this implies  $b \in A_{0}$, in case  $b \neq 0$
(obviously in case $b=0$ we also have $b \in A_{0}$). Hence $K=yA
\cap A_{0}=yA_{0}$, this contradiction with cyclicity   completes
the proof.
\end{proof}

{\bf Remark.}
Let us mention here that for the question on existence of
non-trivial stably free modules it is enough to look only at ideals.
In other words the existence of non-trivial stably free (right)
module is equivalent to the existence of non-trivial stably free
(right) ideal. This follows from the simple observation that if we
have a f.g. projective module, which is non-free, then we also have
a projective non-free ideal. Indeed, let $P$ be a f.g. projective
$R$-module and $R^n=P \oplus Q$. For any $x \in R^n$ we have a
unique decomposition $x=x_P+x_Q$. Consider a submodule of $R^n$ of
the form $(0) \times... \times(0) \times R
 \times(0) \times... \times(0)=R_j \subset R^n$, and define with respect to the above
decomposition submodules $P_j=\{x_P \,| \, x\in R_j \} \subset P$
and $Q_j=\{x_Q \,| \, x\in Q_j \} \subset Q$. Clearly, we have an
isomorphism $R_j = P_j \oplus Q_j$. On the other hand from the
definition of $R_j$ it is clear that $R_j = R$. Moreover we have
$P=P_1 \oplus ... \oplus P_n$ and if $P$ is not free, then one of
$P_j$ is not free. But $P_j$ is a submodule of $R$, i.e. an ideal in
$R$. So we get a projective ideal which is non-free. This remark
shows why it is enough to consider in the theorem only behavior of
ideals under the extension of the base ring.

\medskip

Now we will formulate filtered versions which will be used for the
results about crossed products. Here we restrict ourselves by an
arbitrary well-ordered semigroup.

\begin{lemma}\label{le1f} Let $A=\cup_{\sigma\in \Sigma}U_{\sigma}$    be
a domain, endowed with a filtration by a well-ordered semigroup
$\Sigma$. Then $U_0$ is a $s.c.p.$ -- subalgebra of $A$.
\end{lemma}

\begin{theorem}\label{ft3}
Let $A=\cup_{\sigma\in \Sigma}U_{\sigma}$    be a domain, endowed
with a filtration by a well-ordered semigroup $\Sigma$, and $A$ is
faithfully flat as a left $U_0$-module. Then any stably free
non-free right ideal in $U_0$ can be lifted to $A$.
\end{theorem}

Proofs are  analogous to those of lemmas \ref{le1}, \ref{le2} and
Theorem \ref{gt3}.

\section{Back to crossed products}

Now we can prove the second lemma we need  in the cross product
case.

\begin{lemma}\label{l2} A subalgebra $A_1=A \star U{\cal G}_1$ in $B=A
\star U{\cal G}$, where ${\cal G}_1 \subset {\cal G}$ is a Lie
subalgebra of ${\cal G}$ generated by one (nonzero) element, is a
$s.c.p.$ - subalgebra.

\end{lemma}

\begin{proof} Most essential point in this proof is an appropriate
choice of filtration on $B$. After that we apply lemma \ref{le1f}.
Instead of using a natural filtration on $B$ mentioned at the end of
section \ref{scp}, we suggest the following one. Let
$$B= \bigcup_{\bar i \in
\Sigma_{n-1}} B_{\bar i},
$$

where
$$B_{\bar i}=A[g_1] U_{\bar i},$$

for

 $$U_{\bar i}={\rm Sp}
 \langle g_2^{i_2}...g_n^{i_n  } \, | \, \bar i = (i_2,...,i_n) \in
 \Sigma_{n-1}\rangle_{\k},$$

%  $$\bar i \in
% \Sigma_{(i_2,...,i_n)}=\{(i_2,...,i_n)\,|\, i_k \in \Z_+\}$,

in particular, $B_0 = A[g_1]$ is a polynomial algebra over A on one
variable $g_1$, where ${\cal G}_1$ generated by $g_1$ (as we set in
section \ref{int}).

Note that it is a filtration by well-ordered semigroup
$\Sigma_{n-1}$ and an order on it is degree-lexicographical (the
same as we used for ordering of types in section \ref{scp}, but this
time with respect to $n-1$ variables $g_2,...,g_n$). It is an easy
exercise then to check that this is indeed a filtration.

\end{proof}

Using tools provided by the lemmas \ref{l1} and \ref{l2} we can lift
nontrivial stably free modules from subalgebras of the type $A_1=A
\star U{\cal G}_1=A[g,\delta]$ in a crossed product.

\begin{theorem}\label{t3}
Let $B=A \star U{\cal G}$. Let $K$ be a  nontrivial stably free
right ideal in $A_1=A \star U{\cal G}_1=A[g,\delta]$, for some $g
\in {\cal G}$. Then the induced right ideal $K \otimes_{A_1} B$ in
$B$ is stably free, but not free.
\end{theorem}

\begin{proof}
We use the same filtration as in the previous lemma and apply
theorem \ref{ft3} together with lemmas \ref{l1} and \ref{l2}   for
it.
\end{proof}

The lifting technique could be applied whenever we have a
$s.c.p.$--subalgebra $D$ in $B$, such that ${}_D B$ is a faithfully
flat module. Lemmas \ref{l1} and  \ref{l2} ensure that it is always
the case for the crossed  product algebra $B=A \star U{\cal G}$, if
we choose as a subalgebra $D$ a simple Ore extension $A[g, \delta]$
of $A$.

Now we are in a position to state the result which gives a
sufficient condition of existence of nontrivial stably free modules
over crossed products.

\begin{theorem}
Let $A$ be a noetherian domain, $U{\cal G}$ - the universal
enveloping of Lie algebra ${\cal G}$, and $B=A \star U{\cal G}$ a
crossed product. If there exists an element $g \in {\cal G}$
such that $(r,g+q)$ is a unimodular row in  a subalgebra $A[g,
\delta]$ of $B$, for some $r,s \in A$, $r$ a non-unit, then the
ideal $rB \cap (g+q)B$ is a non-trivial stably free $B$-module.
\end{theorem}

This result shows that nontrivial stably free modules can be
lifted from the Ore extensions of the basic ring $A$, appeared
inside the construction of the crossed product with the universal
enveloping. Obviously these modules not always exist over $A \star
U{\cal G}$.

This we can see already from the example of a simple Ore extension
$A[g, \delta]$, which is also a simplest case of a crossed product.
Let take $A$ to be a commutative local ring with the maximal ideal
$\mu$, it is known \cite{staff} that a nontrivial Ore extension of
$A$ allow stably free non-free ideals if and only if at least one of
the following conditions fails: 1). ${\rm Kdim} A = 1$ or 2).
$\delta(\mu) \subseteq \mu$. Thus situation in the wider class of
crossed products is not so definitive as in group algebras of
solvable groups or in Weyl algebras where nontrivial stably free
non-free modules always exist, so we only can give conditions when
they do.

Let us mention also the following immediate corollary of the
mentioned above fact and theorem \ref{t3}.

\begin{corollary} A crossed product of a local commutative ring $A$
of ${\rm Kdim } A > 1$ with $U \cal G$ for an arbitrary Lie algebra
$G$ always allow stably free non-free module. If  $\,{\rm Kdim } A
\leq 1$ then nontrivial stably free module does exist if  $\,\cal G$
acts in such a way that for some $g \in \cal G$, $g(\cal M)
\not\subset {\cal M }$, where $\cal M$ is a maximal ideal in $A$.
\end{corollary}

\section{Remark on modules of higher ranks}

Here we remind some known results, just to emphasize that  in the
class of crossed products there are obvious examples of non-trivial
stably free modules of higher ranks. They can be obtained  by a
slight modification of arguments for the case of 2-sphere (see
\cite{McR}, 11.2.3).

Namely, let us take
 a (commutative) ring $A={\mathbb
R}[x_1,...,x_n]/ \sum\limits_{i=1}^n x_i^2-1$, for $n \geq 3$. Due
to the nature of these relations the column $\begin{pmatrix}a_1\\ \vdots\\a_n \end{pmatrix}$, with entries $a_i$ ---
images of variables $x_i$ under the natural morphism $\phi: {\mathbb
R}[x_1,...,x_n] \to A$, is unimodular, that is  $Aa_1+...+Aa_n=A$.

 Hence it defines a split monomorphism

$$\alpha: A \hookrightarrow A^n: a \mapsto \begin{pmatrix}a_1\\ \vdots\\a_n \end{pmatrix} \cdot  a$$
with cokernel $P$, so $P \oplus A = A^n$.
Suppose that $P$ is a free $A$-module. This is equivalent to the
fact that the column $\begin{pmatrix}a_1\\ \vdots\\a_n \end{pmatrix}  $ is extendable to an invertible
matrix. That is, there exists $M \in Gl_n({\mathbb R}), M=(\bar r,
\bar c_1,..., \bar c_{n-1})$, where $\bar r, \bar c_1,..., \bar
c_{n-1}$ denote columns of the matrix and $\bar r=\begin{pmatrix}a_1\\ \vdots\\a_n \end{pmatrix}$.
We can construct a continuous tangent vector field on a sphere
${\mathbb S}^{n-1}$, it is provided by the minors $v_r=M_{i2}(r)$ of
matrix $M$, corresponding to the second column. Indeed, the scalar
product $(r,v_r)= {\rm det} (\bar r,\bar r, \bar c_2,..., \bar
c_{n-1})=0$. On the other hand this vector field can't vanish, since
there exist a vector $\bar c_1$, such that $(c_1,v_r)={\rm det}
(r_1,c_1,c_2,...,c_{n-1})={\rm detM} \neq 0$.

The existence of a continuous tangent vector field on a real $n-1$
sphere which vanishes nowhere does contradict, for even $n$, with the
well known theorem on the "brushing of a hedgehog" (or "hairy ball
theorem", see for example \cite{Milnor}).

\section{Remark on non-gradable modules}

Let us mention that for the class of RIT algebras, which form a
special case of crossed products, we can state that nontrivial
stably free modules, we construct here, are also examples of
non-gradable modules.

% This is a consequence of Auslander-regularity
%of RIT algebras.

\bigskip

 {\bf Acknowledgments:}
We would like to thank DFG for the support of the collaborative
research project N30137 and the Phythagoras II project of
Eurocommission.

\end{document}